\documentclass[12pt,reqno]{article}

\usepackage[usenames]{color}
\usepackage{amssymb}
\usepackage{graphicx}
\usepackage{amscd}

\usepackage[colorlinks=true,
linkcolor=webgreen,
filecolor=webbrown,
citecolor=webgreen]{hyperref}

\definecolor{webgreen}{rgb}{0,.5,0}
\definecolor{recrown}{rgb}{.6,0,0}

\usepackage{color}
\usepackage{fullpage}
\usepackage{float}

\usepackage{graphics,amsmath,amssymb}
\usepackage{amsthm}
\usepackage{amsfonts}
\usepackage{latexsym}

\begin{document}


\begin{center}
\vskip 1cm{\LARGE\bf Log-concavity of Lucas sequences of first kind
}
\vskip 1cm
\large
Piero Giacomelli\\
Dipartimento di matematica pura e applicata\\
Universit\`{a} degli studi di Padova\\
Italy\\
\href{mailto:pgiacome@gmail.com}{\tt pgiacome@gmail.com} \\
\end{center}

\vskip .2 in

\begin{abstract}
In these notes we address the study of the log-concave operator acting on Lucas sequences of first kind.
We will find for which initial values a generic Lucas sequence is log-concave, and using this we show when the same sequence is infinite log-concave. The main result will help to fix the log-concavity of some well known recurrent sequences like Fibonacci and Mersenne numbers. Some possible generalization for a complete classification of the log-concave operator applied to general linear recurrent sequences is proposed.
\end{abstract}

\theoremstyle{plain}
\newtheorem{theorem}{Theorem}
\newtheorem{corollary}[theorem]{Corollary}
\newtheorem{lemma}[theorem]{Lemma}
\newtheorem{proposition}[theorem]{Proposition}

\theoremstyle{definition}
\newtheorem{definition}[theorem]{Definition}
\newtheorem{example}[theorem]{Example}
\newtheorem{conjecture}[theorem]{Conjecture}

\theoremstyle{remark}
\newtheorem{remark}[theorem]{Remark}

\def\loperator{log-operator  }
\def\infconcave{$\infty$-log concave }
\def\oneconcave{$1$-fold log-concave }
\newcommand{\im}{\mbox{\upshape Im}}

\newcommand{\proofstart}{\begin{proof}}
\newcommand{\tset}{\mbox{$\cal T$}}
\newcommand{\proofend}{\end{proof}}

\newtheorem{Theorem}{Theorem}[section]
\newtheorem{Proposition}[Theorem]{Proposition}
\newtheorem{Corollary}[Theorem]{Corollary}
\theoremstyle{definition}
\newtheorem{Example}[Theorem]{Example}
\newtheorem{Remark}[Theorem]{Remark}
\newtheorem{Problem}[Theorem]{Problem}

\section{Introduction}
\label{Introduction}
Log-concave sequences arise in many areas of algebra, combinatorics, and geometry as detailed by  the survey article of Brenti \cite{Brenti}.
During the years there have been some studies on the \loperator $\mathcal{L}$ acting on recurrent sequences such as the work of Asai \cite{Asai} on Bell numbers, B\'{o}na \cite{Bona} on sequences counting permutations, Liu \cite{Liu2007453} on combinatorial sequences and McNamara\cite{McNamara20101} with his work on Pascal's triangle. Lucas sequences where first introduced in 1874 by the French mathematician Edouard Lucas, an extensive reference is the book of Koshy \cite{Koshy}.
By definition let $P,Q$ two integer numbers such that $P^2-4Q \geq 0$, then the Lucas sequence of first kind $U_n(P,Q)$ is the recurrent sequence defined by
$U_0 = 0, U_1 = 1, U_2= p,U_n = PU_{n-1} - QU_{n-2}$. As special case for some $P,Q$ the Lucas sequence associated becomes a well known sequence, for example
$L(1,-1,n) = F_n$ where $F_n$ is the Fibonacci sequence. In these notes we study the \loperator on these sequence to address the general problem to find which $P,Q$ integer the corresponding Lucas sequence $U_n(P,Q)$ is log-concave or \infconcave. In section one we will introduce some basic definition and some basic results on \loperator acting on  recurrent sequences.
Section two will show a general result on how to solve the log-concavity problem on a generic Lucas sequence of first kind.
Last section will propose a generalization of the methods used on Lucas sequence to  generic linear recurrent sequences.

\section{Basic definition}
\label{Basic definition}
We now remark some definitions of the \loperator. We refer to the notation to McNamara \cite{McNamara20101}. Let us start with
\begin{definition}
Let $(a_n)_{n\in\mathbb{N}}$ a real sequence we define the \loperator  as a function $\mathcal{L}: \mathbb{R}\rightarrow\mathbb{R}$ such that $b_n = \mathcal{L}(a_n) = a_n^2 -a_{n-1}a_{n+1}$. If $b_n \geq 0 $ for all $n \in \mathbb{N}$ then the sequence $(a_n)_{n \in \mathbb{N}}$ is said to be log-concave.
\end{definition}
Considering that log-concavity can deals with negative indexes, by convention we will extend a sequence $(a_n)_{n\in\mathbb{N}}$ to a sequence $(a_n)_{n\in\mathbb{Z}}$ where by definition $a_n = 0$ if $n <0$. In the same way if the sequence is finite so $n \leq m, m \in \mathbb{N}$ then all other indexes $n > m$ will be zero.

\begin{definition}
A real sequence $(a_n)_{n\in\mathbb{N}}$ is said to be $i-fold$ log-concave for $i\in\mathbb{N}, i \geq 1$ if  $\mathcal{L}^i(a_n)$ is a nonnegative sequence. Where $\mathcal{L}^i(a_n)$ is the \loperator applied to a sequence $(a_n)_{n\in\mathbb{N}}$ $i$-times so $\mathcal{L}^i = \mathcal{L}\circ\mathcal{L}\circ\dots\circ\mathcal{L}$.
\end{definition}
Using McNamara \cite{McNamara20101} notation:
\begin{definition}
a real sequence $(a_n)_{n\in\mathbb{N}}$ is said to be \infconcave if  $\mathcal{L}^i(a_n)$ is a nonnegative sequence for all $i\in \mathbb{N}, i \geq 1 $.
\end{definition}
So log-concavity in the ordinary sense is $1$-fold log-concavity.
To study log-concavity on Lucas sequences, we need some preliminary results, like the following:
\begin{lemma}
\label{lem:zeroinfinity}
let $(a_n)_{n \in \mathbb{N}}$ a  sequence where $a_n=k$ for all $n\in\mathbb{N}$ and $k$ is a real number, then $(a_n)_{n \in \mathbb{N}}$ is \infconcave.
\end{lemma}
\proofstart
It is easy to check that
\begin{equation*}
b_n= \mathcal{L}(a_n) = \mathcal{L}(a_n) = k^2 - (k\cdot k) = 0
\end{equation*}
for all $n \in \mathbb{N}$. It is also clear that the all zeros sequence $b_n$ is invariant by the \loperator that is $\mathcal{L}(b_n) = b_n$. Being $b_n \geq 0$ that means that also $\mathcal{L}(b_n) \geq 0$ so the all zeros sequence is \infconcave.
\proofend

In the same way it is also easy to check that

\begin{lemma}
\label{recurrentinfinity}
let $(a_n)_{n \in \mathbb{N}}$ a  sequence where for all $n\in\mathbb{N}$ $a_n=kb^n$ where $k,b \in \mathbb{R}, k \neq 0, b \neq 0$  then $(a_n)_{n \in \mathbb{N}}$ is \infconcave.
\end{lemma}
\proofstart
By direct check
\begin{equation*}
\mathcal{L}(a_n) =  (a_n)^2 - a_{n-1} a_{n+1} = k^2b^{2n} - k^2b^{n-1+n+1} = k^2b^{2n} - k^2b^{2n} = 0
\end{equation*}
for all $n \in \mathbb{N}$. So $a_n$ is \oneconcave and the result sequence is the all zeros sequence than
considering lemma \ref{lem:zeroinfinity} then the sequence $a_n$ is also \infconcave.
\proofend

In next section we will detail our analysis of the \loperator to the Lucas sequence.

\section{Log-operator and Lucas sequences}
In these section we address the study of log-concavity, of a Lucas sequence of first kind
Let start with the Lucas sequence definition:
\begin{definition}
Let $(P,Q) \in \mathbb{Z}\times\mathbb{Z}$ two non-zero integer such that $P^2-4Q \geq 0$ and let $n\in\mathbb{N}$ an index. A Lucas sequence $U_n(P,Q) $ of first kind is a recurrent sequence defined as follows:
\begin{align*}
U_0 &= 0 \\
U_1 &= 1 \\
U_n &= PU_{n-1}-QU_{n-2}.
\end{align*}
\end{definition}	
Choosing the correct $P,Q$ it is possible to obtain some well known sequences for example:
\begin{itemize}
	\item{If $P=1,Q=-1$ then the Lucas sequence $U_n(1,-1)=F_n$ where $F_n$ is the Fibonacci sequence. }
\item{If $P=2,Q=-1$ then the Lucas sequence $U_n(2,-1)$ is the sequence of Pell numbers. }
\item{If $P=1,Q=-2$ then the Lucas sequence $U_n(1,2)$ is the sequence of Jacobsthal numbers. }
\item{If $P=3,Q=2$ then the Lucas sequence $U_n(3,2)$ is the sequence of Mersenne numbers. }
\end{itemize}

The main result of this section will prove for which initial $P,Q$ the resulting Lucas sequence is \infconcave. Let us start by showing that in general if we choose a generic couple $P,Q$ the Lucas sequence $U_n(P,Q)$ is not \oneconcave.

We use the following proposition

\begin{proposition}
\label{fibonacci}
The Fibonacci sequence $F_n$ is not \oneconcave.
\end{proposition}
\proofstart
Considering the \loperator applied to $F_n$ we have
\begin{equation*}
b_n = \mathcal{L}(F_n) = F_n^2-F_{n-1}F_{n+1};
\end{equation*}
now by the Cassini's identity
\begin{equation}
F_{n-1}F_{n+1}- F_n^2 = (-1)^n
\end{equation}
we obtain
\begin{equation*}
  F_n^2-F_{n-1}F_{n+1} = (-1)\cdot(-1)^n = (-1)^{n+1}.
\end{equation*}
So
\begin{equation*}
  \mathcal{L}(F_n) = (-1)^{n+1}.
\end{equation*}

thus $F(n)$ is not \oneconcave.
If we applied the $\mathcal{L}$ operator to the sequence $b_n$ and we calculate $\mathcal{L}^2(F(n))=\mathcal{L}(\mathcal{L}(b_n)$ we obtain
\begin{equation*}
\mathcal{L}^2(F(n)) = ((-1)^{n+1})^2-(-1)^{n+2}\cdot(-1)^{n} = ((-1)^{n+1})^2-(-1)^{2n+2} = 1 - 1 = 0
\end{equation*}
so after applying the \loperator more than once we obtain a sequence that is log-concave.
\proofend

We will now fix for what initial parameter $P,Q$ the generate Lucas sequence $U_n(P,Q)$ is a $1$-fold log-concave Lucas sequence, and in these cases where for what $P,Q$ the Lucas sequence becomes \infconcave. Instead of trying to apply directly the \loperator to the generic expression of the Lucas sequence $U_n(P,Q)$, we will use a more treatable expression for $U_n(P,Q)$. To do this, we first need to recall \cite{Binet} that:
\begin{remark}
\label{mainremark}
let $U_n(P,Q)$ a Lucas sequence of first kind, than the characteristic equation of the recurrence relation is
\begin{equation}
x^2-Px+Q=0
\end{equation}
that has discriminant $D=P^2-4Q$. If the discriminant is positive so $D \geq 0$ then the roots of the characteristic equation
are
\begin{equation}
a=\frac{P+\sqrt{D}}{2}, \ b=\frac{P-\sqrt{D}}{2}
\end{equation}
and so if $D \geq 0$ it is possible to rewrite $U_n(P,Q)$ in the following way
\begin{equation}
U_n(P,Q) = \frac{a^n-b^n}{a-b} = \frac{a^n-b^n}{\sqrt{D}}.
\end{equation}

\end{remark}
Armed with this expression for Lucas sequence, we will divide our study in two main cases let us start with the simpler one.
\begin{proposition}\label{prop:zero}
Let $U_n(P,Q)$ a Lucas sequence where $P,Q$ are two integer and the discriminant $D$ of the characteristic equation associated with $U_n(P,Q)$ is zero then the  Lucas sequence associated is \oneconcave.
\end{proposition}
\proofstart
It is easy to see that if $D = 0$ then $P^2-4Q=0$ and so there exists and integer $S$ such that $P=2S$ and $Q=S^2$. Using this fact the Lucas sequence associated can be rewritten in the form
\begin{equation}
U_n = n S^{n-1}.
\end{equation}
So now, applying the $\mathcal{L}$ operator, we see that
\begin{align*}
\mathcal{L}(U_n) &= \\
\mathcal{L}(nS^{n-1}) &= (nS^{n-1})^2-[(n-1)S^{n-2}\cdot(n+1)S^n]   \\
											&= (n^2)S^{2n-2}-(n^2-1)S^{n-2+n}   \\
											&= (n^2)S^{2n-2}-(n^2-1)S^{2n-2}   \\
											&= (n^2-n^2+1)S^{2n-2}   \\
											&= (S^{n-1})^2
\end{align*}
and so $\mathcal{L}(U_n) \geq 0$ for all $S \in \mathbb{Z}$. This prove that $U_n$ is \oneconcave.
\proofend
From proposition \ref{prop:zero} we have also the following corollary
\begin{corollary}
Let $U_n(P,Q)$ a Lucas sequence where $P,Q$ are two integer and there exist an $S \in \mathbb{Z}$ such that $P=2S$ and $Q=S^2$  then the
Lucas sequence associated is \infconcave.
\end{corollary}
\proofstart
We have seen that under the hypothesis $\mathcal{L}(U_n) = (S^{n-1})^2 = (S^2)^{n-1}$.
By changing  the index we have that the original sequence become a  sequence of the form $b_k = {S^k}$
where $k \in \mathbb{Z}, k = 2n-2, k \geq -2 $.
Considering that for negative indexes
$b_k = 0$ we have that by lemma \ref{recurrentinfinity} the sequence $b_k$  is \infconcave and so $U_n$.
\proofend

Let now consider the general case

If $D=P^2-4Q > 0$ by remark \ref{mainremark} it is possible to rewrite $U_n(P,Q)$ in the following way
\begin{equation}
U_n(P,Q) = \frac{a^n-b^n}{a-b} = \frac{a^n-b^n}{\sqrt{D}}
\end{equation}
where
\begin{equation}
a=\frac{P+\sqrt{D}}{2}, \ b=\frac{P-\sqrt{D}}{2}
\end{equation}

we notice that, using direct calculation we have
\begin{align*}
\mathcal{L}(U_n) &= \\
\mathcal{L}\biggl( \frac{a^n-b^n}{\sqrt{D}} \biggl) &= \biggl(\frac{a^n-b^n}{\sqrt{D}}\biggl)^2 - \biggl[\frac{a^{n-1}-b^{n-1}}{\sqrt{D}}\cdot\frac{a^{n+1}-b^{n+1}}{\sqrt{D}}\biggl]   \\
&= \biggl(\frac{a^{2n}-2a^nb^n+b^{2n}}{D}\biggl)
-
\biggl(\frac{a^{n-1+n+1} -a^{n-1}b^{n+1} -a^{n+1}b^{n-1} +b^{n+1+n-1} }{D}\biggl)
\\
&= \frac{a^{2n}-2a^nb^n+b^{2n} - a^{2n} +  a^{n-1}b^{n+1} +a^{n+1}b^{n-1} -b^{2n}  }{D}
\\
&= \frac{a^{n+1}b^{n-1}   -2a^nb^n +  a^{n-1}b^{n+1}   }{D} \\
&= \frac{ a^{n-1}b^{n-1} ( a^2   -2ab +  b^2  ) }{D} \\
&= \frac{ (ab)^{n-1} ( a^2   -2ab +  b^2 )  }{D} \\
&= \frac{ (ab)^{n-1} ( a-b )^2  }{D} \\
\end{align*}
now then by definition
\begin{equation}
ab = \frac{P+\sqrt{D}}{2}\cdot \frac{P-\sqrt{D}}{2}  = \frac{1}{4} (P^2-D) = \frac{1}{4} (P^2-P^2+4Q) = Q
\end{equation}
and
\begin{equation}
a-b = \frac{P+\sqrt{D}}{2}  -  \frac{P-\sqrt{D}}{2}  = \frac{2P}{2} = P.
\end{equation}
So finally we have

\begin{equation}
\mathcal{L}(U_n) = \frac{Q^{n-1}P^2}{D}
\end{equation}
So $\mathcal{L}(U_n) \geq 0$ if $Q \geq 0$. Combining this with the assumption that $P^2-4Q \geq 0$ we have that $U_n(P,Q)$ is 	 \oneconcave if
\begin{equation*}
\begin{cases}
Q \geq 0 \\
P^2-4Q > 0
\end{cases}
\end{equation*}
that gives the following set of solutions
$Q \geq 0 \wedge P > 2\sqrt{Q}$ or $Q \geq 0 \wedge P < - 2\sqrt{Q}$.

We can summarize the result in the following
\begin{theorem}\label{thm:main}
Let $P,Q$ two integer such that $Q \geq 0 \wedge P > 2\sqrt{Q}$ or $Q \geq 0 \wedge P < - 2\sqrt{Q}$, then the associated Lucas sequence $U_n(P,Q)$ is \oneconcave.
\end{theorem}
using theorem \ref{thm:main} and the lemma \ref{recurrentinfinity}, it is easy to check that
\begin{corollary}\label{cor:final}
Let $P,Q$ two integer such that $Q \geq 0 \wedge P > 2\sqrt{Q}$ or $Q \geq 0 \wedge P < - 2\sqrt{Q}$. Then the Lucas sequence $U_n(P,Q)$ is \infconcave.
\end{corollary}
\proofstart
Under the hypothesis we have that
\begin{equation*}
b_n =\mathcal{L}(U_n) = \frac{Q^{n-1}P^2}{D}
\end{equation*}
that is a sequence of the form $k_na^n$ and  by lemma \ref{recurrentinfinity} $U_n(P,Q)$ is \infconcave.
\proofend
At the end using the corollary \ref{cor:final} we can check that:
\begin{itemize}
\item{$U_n(1,-1)$ is the Fibonacci sequence that is not \oneconcave and so neither \infconcave.}
\item{$U_n(2,-1)$ is the sequence of Pell numbers that is not \oneconcave and so neither \infconcave . }
\item{$U_n(1,-2)$ is the sequence of Jacobsthal numbers that is not \oneconcave and so neither \infconcave. }
\item{$U_n(3,2)$ is the sequence of Mersenne numbers that is \infconcave. }
\end{itemize}

\section{Conclusion}
In these notes we have studied the \loperator applied to a generic Lucas sequence of first kind $U_n$. We have shown that for initial parameter $Q \geq 0, P \geq 2Q $ or $Q \geq 0, P \leq 2Q $, the associate Lucas sequence of first kind is \infconcave. As result we find that Fibonacci, Pell and Jacobsthal sequences are not \infconcave but the Mersenne numbers sequence is \infconcave.
There is a natural question that arise from these results. As shown the key fact, that a sequence is recurrent, allow the sequence to be expressed in a more treatable way before applying the \loperator. It would be interesting giving a generic linear recurrent sequence that satisfy a generic characteristics equation of order $k$, to find sufficient condition on the coefficient of the equation to be sure that the sequence is \oneconcave and after this fix which conditions leads to a \infconcave sequence. Formalizing a little, giving a recurrent sequence define as $a_n= k_1a_{n-1} + k_2a_{n-2}+\dots k_ma_{n-m}$ that has a characteristic equation  $a_n - k_1a_{n-1} - k_2a_{n-2}-\dots-k_ma_{n-m} = 0 $ is there is a sufficient condition on the $k_1,k_2,\dots,k_m$ integer coefficient such that $a_n$ is \oneconcave and \infconcave. This question would be subject of further study.

\bigskip
\hrule
\bigskip


\begin{thebibliography}{99}


\bibitem[1]{Brenti}
F. Brenti, Log-concave and Unimodal sequences in Algebra, Combinatorics, and Geometry: an update, \textit{Contemporary Mathematics}.\textbf{178} (1994),71--89.


\bibitem[2]{Asai}
N. Asai, I. Kubo and H. Kuo, Bell numbers, log-concavity, and log-convexity, \textit{Acta Applicandae Mathematicae}. \textbf{63} (2000), 79--87.

\bibitem[3]{Bona}
M. B\'{o}na, A Combinatorial Proof of the Log-Concavity of the Numbers of Permutations with $k$ Runs, \textit{J. Combin. Theory Ser. A}. \textbf{90} (2000), 293--303.

\bibitem[4]{Liu2007453}
L. L. Liu and Y. Wang, On the log-convexity of combinatorial sequences, \textit{Advances in Applied Mathematics} \textbf{39} (2007), 453--476.

\bibitem[5]{McNamara20101}
P. R. W. McNamara and B. E. Sagan,Infinite log-concavity: Developments and conjectures, \textit{Advances in Applied Mathematics} \textbf{44} (2010), 1--15.

\bibitem[6]{Koshy}
T. Koshy, \textit{Fibonacci and Lucas Numbers with Applications}, Wiley, John \& Sons, Incorporated , 2001

\bibitem[7]{Binet}
H. W. Austin and J. W. Austin, Binet Formulas for Recursive Integer Sequences, \textit{Journal of Mathematical sciences and Mathematical education} \textbf{4} (2009).


\end{thebibliography}
\end{document}